\documentclass[12pt]{article}

\scrollmode

\usepackage{amsmath}
\usepackage{amsfonts}

    \newtheorem{theorem}{Theorem}[section]
    \newtheorem{lemma}[theorem]{Lemma}

\newenvironment{proof}[1][Proof]{\begin{trivlist}
\item[\hskip \labelsep {\bfseries #1}]}{\end{trivlist}}
\newcommand{\qed}{\nobreak \ifvmode \relax \else
      \ifdim\lastskip<1.5em \hskip-\lastskip
      \hskip1.5em plus0em minus0.5em \fi \nobreak
      \vrule height0.75em width0.5em depth0.25em\fi}

\def\eps{\varepsilon}
\def\qed{\hfill\rule{.2cm}{.2cm}}
\def\P{{\mathbb P}}
\def\esp{{\mathbb E}}
\def\Z{{\mathbb Z}}
\def\N{{\mathbb N}}

\def\1{{\mathbf 1}}

\def\a{\alpha}
  
\def\l{\lambda}

\def\g{\gamma}

\def\x{\Gamma}
\def\G{\Gamma}

\def\d{\delta}

\def\cp{{\cal P}}

\def\t{{\mathfrak T}}

\def\p{{\cal P}}

\def\RR{{\cal R}}

\def\D{{\cal D}}

\def\x1{X_1^{(1)}}

\def\tx0{\tilde X^0}

\def\={&=&}
\def\+{&+&}

\def\ja{\xrightarrow{J_1}}
\def\jpa{\xrightarrow{J_1,P}}
\def\jppa{\xrightarrow{J_1,P}}
\def\joa{\xrightarrow{J_1,P_1}}
\def\jta{\xrightarrow{J_1,P_2}}
\def\jtha{\xrightarrow{J_1,P_3}}

\def\jfia{\xrightarrow{J_1,P_5}}

\newtheorem{theo}{Theorem}%[section]

\newtheorem{cor}[theo]{Corollary}
\newtheorem{rmk}[theo]{Remark}

\def\beq{\begin{equation}}
\def\eeq{\end{equation}}
\newcommand{\bei}{\begin{itemize}}
\newcommand{\eei}{\end{itemize}}
\newcommand{\ben}{\begin{enumerate}}
\newcommand{\een}{\end{enumerate}}
\newcommand{\beqn}{\begin{eqnarray}}
\newcommand{\beqnn}{\begin{eqnarray*}}
\newcommand{\eeqn}{\end{eqnarray}}
\newcommand{\eeqnn}{\end{eqnarray*}}
\newcommand{\brm}{\begin{rmk}}
\newcommand{\erm}{\end{rmk}}

\usepackage{amsmath}
\usepackage{amsfonts}

\begin{document}

\title{Scaling limits and aging for asymmetric trap models on the complete graph and $K$ processes}
%\date{}
\author{
S.~C.~Bezerra \footnote{CCEN-UFPe, Cidade Universitária, 50740-540
Recife PE, Brazil, lrenato@ime.usp.br} \thanks{Supported by 
FAPESP fellowship 2007/03517-3}
\and
L.~R.~G.~Fontes \footnote{IME-USP, Rua do Mat\~ao 1010, 05508-090
S\~ao Paulo SP,  Brazil, lrenato@ime.usp.br} \thanks{Partially
supported by CNPq grant 305760/2010-6, and FAPESP grant 2009/52379-8}
\and
R.~J.~Gava \footnote{IME-USP, Rua do Mat\~ao 1010, 05508-090 S\~ao Paulo SP,  Brazil, gavamat@yahoo.com.br} 
\thanks{Supported by FAPESP fellowship 2008/00999-0}
\and
V.~Gayrard\footnote{CMI, 39 rue Joliot Curie, 13453 Marseille, France, veronique@gayrard.net}
\mbox{ }\thanks{Partially supported by FAPESP grants 2007/59096-6 and 2009/51609-0}
\and
P.~Mathieu\footnote{CMI, 39 rue Joliot Curie, 13453 Marseille, France,
pierre.mathieu@cmi.univ-mrs.fr}  
}
\date{\today}  
 
\maketitle

\pagestyle{myheadings}
\markright{Scaling limits for asymmetric trap models and $K$ processes}

 \begin{abstract} We obtain scaling limit results for asymmetric trap models and their infinite volume counterparts, namely asymmetric
$K$ processes. Aging results for the latter processes are derived therefrom.
 \end{abstract}


\smallskip


\section{Introduction}

The long time behavior of trap models and related processes with disordered parameters has been the theme of several papers in the recent literature.
From the inaugurating work of Bouchaud~\cite{B}, where the case of the complete graph was shown to exhibit {\em aging}, the same as well as other
cases were analysed. The model on the complete graph was further studied in~\cite{BF} and~\cite{FM1}, with different points of view, and considering
distinct time scales. And more recently,~\cite{G} took up the asymmetric case, which is also the model we study here.

The trap model in the complete graph is sometimes also called REM-like trap model, due to its resemblance to a dynamics for
the Random Energy Model (REM~\cite{D}). Such a dynamics for the REM, on the hypercube rather than the complete graph, was studied
in~\cite{ABG1,ABG2}, where aging results comparable to the ones of Bouchaud were derived.  See also~\cite{C1,FL,G1}.
Trap and trap-like models associated to correlated energy (mean field) spin glasses have been the object of more recent work: 
a dynamics for the $p$-spin model was studied in~\cite{ABC,BG}, and results on the GREM-like trap model were obtained in~\cite{FGG}.

Trap models on $\Z^d$ have also attracted a lot of interest, in connection with aging as well as with {\em localization}; 
see~\cite{FIN,AC1,AC2,ACM,FM2} -- results on the asymmetric case were obtained recently in~\cite{BC,C2,M}. Analyses on tori were
performed in~\cite{AC,JLT}.

In this paper, we revisit the trap model in the complete graph, described briefly below in this introduction, and in full in Section~\ref{sec:tck}. 
Our goal is twofold: 
\begin{enumerate}
 \item to propose a representation of the model -- in terms of trap depth, rather than location -- for which scaling limits can be derived in a
unified manner in different scaling regimes;
 \item and to introduce the infinite volume processes which result from these scaling limits, in particular the asymmetric $K$ process.
\end{enumerate}

Let us now briefly describe the asymmetric trap model in the complete graph with $n$ vertices. This is a 
continuous time Markov chain on the vertices of that graph, whose mean jump time at site $x$ is given by $\tau_x^{1-a}$, where $a\in[0,1]$ is 
an asymmetry parameter, 
and whose transition probability from any site $x$ to any site $y$ is proportional to $\tau^a_y$, where $\{\tau_x\}$ are iid positive random
variables in the domain of attraction of an $\a$-stable law. The random variable $\tau_x^{1-a}$ may be interpreted as the depth of the trap 
at site $x$. 
One readily checks that this dyanmics is reversible with respect to the measure whose weights are given by $\{\tau_x\}$.
The case $a=0$ is that of the {\em symmetric} model. We call the general case where $a\in[0,1]$ the {\em asymmetric} model.
Let $Y_n(t)$ denote the site visited at time $t$.

This paper is more immediately related to~\cite{FM1} and~\cite{G}, so let us briefly outline our results here against the background of the
ones of those papers. 
Whereas in the former reference a scaling limit was derived %as volume goes to infinity 
for the symmetric model at times of the order of the deepest trap in the landscape, and then aging results were derived for a class of 
two-time correlation functions of the limit model at vanishing times, in here we present similar limit results for the asymmetric model.
Rather than looking at $Y_n(t)$ however, we consider $Z_n(t)=\tau^{1-a}_{Y_n(t)}$, the depth of the currently visited trap. 
As explained below, this is a convenient representation for taking scaling limits, not only at times of the order of the deepest trap in the
landscape, which we do here using this representation (see Theorem~\ref{teo:es}), obtaining a limiting process which we denote by $Z$, 
but at shorter time scales as well.
We call $Z$ the {\em asymmetric $K$ process}, in allusion to the $K$ process introduced in~\cite{FM1}.
We further derive a scaling limit result for $Z$ at vanishing times (see Theorem~\ref{teo:asz}), obtaining a limiting process $\hat Z$ 
which is self similar of index 1. 
The latter fact may be interpreted as a fuller aging result for $Z$, involving the dynamics itself, not only a class of correlation 
functions thereof. Other scaling regimes of $Z_n$ may be analysed with the same approach, with similar results.

Scaling limits of asymmetric trap models in the complete graph are also the main theme of~\cite{G}. In that work scaling limits of the {\em
clock process} are derived in several scaling regimes (essentially all of them: from ``order 1'', where the volume limit is taken first, and then
the time limit, to the scale where the model is virtually at equilibrium, including scales in between,
in particular the ones treated here); occurrence of aging and other dynamical phenomena are discussed for each regime. 

One reason to consider a representation like $Z_n$, as we do here, rather than the clock process, is that, besides the information on the jump times
given by the latter process, $Z_n$ provides also location information, absent in that process. For, say, correlation functions which depend only on
jump times (like the $\Pi$ functions discussed on Subsection~\ref{ssec:agezz} below; see~(\ref{eq:ag2a}-\ref{eq:ag2b})), the clock process is enough.
But other ones require location information, and in those cases the clock process is no longer enough on its own. 
We discuss two such examples in Subsection~\ref{ssec:agezz} below.

$Z_n$ and $Z$, as well as their rescaled versions, and $\hat Z$ also, can be described as functions of two related subordinators, 
the second being obtained as the integral of an independent iid family of mean 1 exponentials with respect to the first one. 
Once we obtain the limit of the first subordinator in a given scaling regime, a continuity property of the above 
mentioned function implies a limit result for the original process.
Section~\ref{sec:lm} below is devoted to establishing that continuity property (see Lemma~\ref{lemma}) in a somewhat abstract setting, which may
turn out to be the setting of similar processes of interest. 

In Section~\ref{sec:tck} we describe our trap models and $K$ processes in more detail and then, applying the auxiliary result of
Section~\ref{sec:lm}, 
we derive scaling limit results for them, as anticipated above, the one for the trap model in Subsection~\ref{ssec:erg},  
and the one for the $K$ process in Subsection~\ref{ssec:agez}. In the closing Subsection~\ref{ssec:agezz} we discuss the derivation of aging results
for three particular two-time correlation functions of $Z$ as a corollary to Theorem~\ref{teo:asz}.

\section{A continuity lemma about a class of trajectories in $D$}
\label{sec:lm}
\setcounter{equation}{0}

Let $D$ be the space of c\`adl\`ag real trajectories on $\mathbb{R}^+=[0,\infty)$ equipped with the $J_1$ Skorohod metric
(see e.g.~\cite{EK} Chapter 3, Section 5). Let $\mathbb{N}^*=\{1,2,\ldots\}$ denote the positive integers.

Let $S,\,S^\eps,\,\eps>0,$ be nonnegative nondecreasing jump functions in $D$, i.e., suppose that 
there exist (countable) subsets $A^\eps=\{x_i^\eps,\, i\in \N^* \}$ and $A=\{x_i,\, i\in \N^* \}$ of $\mathbb{R}^+$ 
and positive number sequences 
$\{\gamma^\eps_{x_i^\eps},\, i\in \N^* \}$ and $\{\gamma_{x_i},\, i\in \N^* \}$ such that
\begin{equation}
\label{eq:s} 
S_r^\eps=\sum_{i\colon x_i^\eps\in[0,r]}\gamma^\eps_{x_i^\eps}<\infty,\quad S_r=\sum_{i\colon x_i\in[0,r]}\gamma_{x_i}<\infty,\quad r\geq0.
\end{equation}

Consider $\{T_i,\,i\in \N^*\}$, a family of i.i.d.~exponential random variables of mean 1 and let

\begin{equation}
\label{eq:g} 
\Gamma_r^\eps=\sum_{i\colon x_i^\eps\in[0,r]}\gamma^\eps_{x_i^\eps}T_i,\quad \Gamma_r=\sum_{i\colon x_i\in[0,r]}\gamma_{x_i}T_i,\quad r\geq0,
\end{equation}
%and 
\begin{equation}
\label{eq:ze} 
Z_t^\eps=\begin{cases}
              \gamma^\eps_{x^\eps_{i_0}},&\mbox{  if }t\in[\Gamma^\eps_{x^{\eps}_{i_0}-},\Gamma^\eps_{x^\eps_{i_0}})\mbox{ for some }i_0,\\
\mbox{}\,\,\,0,&\mbox{  if }t\notin[\Gamma^\eps_{x^{\eps}_{i}-},\Gamma^\eps_{x^\eps_{i}})\mbox{ for any }i,
             \end{cases}
\end{equation}
and
\begin{equation}
\label{eq:z} 
Z_t=\begin{cases}
              \gamma_{x_{i_0}},&\mbox{  if }t\in[\Gamma_{x_{i_0}-},\Gamma_{x_{i_0}})\mbox{ for some }i_0,\\
\mbox{}\,\,\,0,&\mbox{  if }t\notin[\Gamma_{x_{i}-},\Gamma_{x_{i}})\mbox{ for any }i.
             \end{cases}
\end{equation}

Below, we will use the symbol $\ja$ to denote (strong) convergence on $(D,J_1)$, while $\jpa$ will denote weak convergence on $(D,J_1)$
with respect to a given probability measure $P$.

Let $P$ denote the probability measure on $(D,J_1)$ induced by the distribution of $\{T_i,\,i\in \N^*\}$.% via the above constructions.

\begin{lemma} 
\label{lemma}
Let $S^\eps,\,S,\,Z^\eps,\,Z$ be as above.
As $\eps\to0$, if $S^\eps\ja S$, then $Z^\eps\displaystyle\jppa Z$.
\end{lemma}

\begin{rmk}
 \label{rmk:int}
From~(\ref{eq:s}-\ref{eq:z}), we see that $Z=\Xi(S,\{T_i,\,i\in \N^*\})$ and $Z^\eps=\Xi(S^\eps,\{T_i,\,i\in \N^*\})$, where $\Xi$ is the composition
underlying the above definitions. Lemma~\ref{lemma} then establishes a continuity property of the distribution of $\Xi$ in its first argument.
\end{rmk}

\noindent{\bf Proof of Lemma~\ref{lemma}}

We will assume that there exists $R'\in\mathbb{R}^+$ such that $|A\cap[0,R']|=\infty$. Other cases may be argued similarly,
with simpler arguments.

Let $\Gamma^{-1}$ be the (right continuous) inverse of $\Gamma$. Let us fix $T>0$. Then one readily checks that, given $\d>0$, 
there exists $R\notin A$, $R\geq R'$, such that
\begin{equation}
\label{eq:d1}
\P(\Gamma^{-1}(T)\geq R)\leq\d.
\end{equation}

Given $\eta>0$, we may choose $\d'>0$ be such that 
\begin{equation}
\label{eq:d2}
S_{R+\d'}-S_R<\eta.
\end{equation}

Let us now enumerate 
$A\cap[0,R]=\{x_1,x_2,\ldots\}$ such that $\gamma_{x_1}\geq\gamma_{x_2}\geq\ldots$.
From the hypothesis, there exists $m=m(\eps)$, with $m\to\infty$ as $\eps\to0$, and an enumeration of
$A^\eps\cap[0,R]=\{x^\eps_1,x^\eps_2,\ldots\}$ such that as $\eps\to0$
\begin{equation}
\label{eq:c1}
\left(\sup_{1\leq i\leq m}|x^\eps_i-x_i|\right)\vee \left(m\sup_{1\leq i\leq m}|\gamma^\eps_{x^\eps_i}-\gamma_{x_i}|\right)\to0.
 \end{equation}
It follows from this and the hypothesis that, given $\eta>0$, for all small enough $\eps$ and $1\leq k\leq m$
\begin{equation}
\label{eq:c1a}
\sum_{i>k}\gamma^\eps_{x^\eps_i}=S^\eps_{R}-\sum_{i=1}^k\gamma^\eps_{x_i^\eps}
\leq S_{R+\d'}-\sum_{i=1}^k\gamma_{x_i}+\eta=S_{R+\d'}-S_R+\sum_{i>k}\gamma_{x_i}+\eta\leq\sum_{i>k}\gamma_{x_i}+2\eta.
 \end{equation}

We now recall that in the $J_1$ topology, functions are close if they are uniformly close inside arbitrary bounded intervals, 
after allowing small time distortions (for details see e.g.~\cite{EK} Chapter 3, Section 5). 

Now, given $k\geq1$ arbitrary but fixed, independent of $\eps$, let $\{\bar x_1,\ldots,\bar x_k\}$ be an enumeration of $\{x_1,\ldots,x_k\}$
such that $\{\bar x_1<\ldots<\bar x_k\}$.  This leads to an enumeration $\{\bar x_1^\eps,\ldots,\bar x_k^\eps\}$ of $\{x_1^\eps,\ldots,x_k^\eps\}$
such that for $1\leq i\leq k$
\begin{equation}
\label{eq:c1b}
\bar x^\eps_i\to\bar x_i\mbox{ and }\gamma^\eps_{\bar x^\eps_i}\to\gamma_{\bar x_i}
 \end{equation}
(see paragraph of~(\ref{eq:c1}) above). At this point we relabel $\{T_i\}$ so that $T_1,\ldots,T_k$ are attached to $\bar x_1<\ldots<\bar x_k$
and commonly to $x_1^\eps,\ldots,x_k^\eps$, respectively, which does not change distributions. Let $Z^{(k)}$ and $Z^{(k,\eps)}$ be the
respective versions
of $Z$ and $Z^\eps$ with the relabeled $\{T_i\}$.

Let us now take a family of temporal distortions $(\lambda^\eps)=(\lambda_{k}^{\eps})$ as follows.
For $1\leq i\leq k$, we consider the time intervals 
$I_i=[t_i^-,t_i]$, where $t_i=\Gamma_{\bar x_i}$ and $t_i^-=\Gamma_{\bar x_i-}$, and $[t_i^{\eps-},t_i^{\eps}]$, where 
$t_i^{\eps}=\Gamma^\eps_{\bar x^\eps_i}$ and $t_i^{\eps-}=\Gamma^\eps_{\bar x^\eps_i-}$,
 % Proposition 5.3 page 119 in \cite{EK}.
and then define
\begin{equation}%{align}
\label{eq:lam}
\lambda^\eps(s)=
\begin{cases}
%\mbox{}\hspace{1cm}
\frac{t_{1}^{\eps -}}{t_{1}^-}\,s,&\mbox{ if }0\leq s\leq t_{1}^-,\\
\frac{t_i^{\eps}-t_i^{\eps -}}{t_i-t_i^-}(s-t_i^-)+t_i^{\eps-},&\mbox{ if }t_i^-\leq s\leq t_i,\\
%\mbox{}\,\,
\frac{t_{i+1}^{\eps -}-t_i^{\eps}}{t_{i+1}^--t_i}(s-t_i)+t_i^{\eps},&\mbox{ if }t_i\leq s\leq t_{i+1}^-,\\
%\mbox{}\hspace{.2cm}
(s-t_{k+1}^-)+t_{k+1}^{\eps-},&\mbox{ if }s\geq t_{k+1}^-,
\end{cases}
\end{equation}%{align}
where $t_{k+1}^-:=\Gamma_R$, $t_{k+1}^{\eps-}:=\Gamma^\eps_R$.

At this point, we have two tasks: the first one is to control the slopes of the functions $\lambda^\eps$ and the second one is to control the 
$\sup$ norm of the difference $Z^{(k,\eps)}_{\lambda^\eps(t)} -Z^{(k)}_t$. 

We start by the second task. Let $\mathcal{M}=\cup_{i=1}^kI_i$. If $t\in\mathcal{M}$, then 
\begin{equation}
\label{eq:c2}
|Z^{(k)}_t-Z^{(k,\eps)}_{\lambda^\eps(t)}|\leq\max_{1\leq i\leq k}|\gamma_{\bar x_i}-\gamma^{\eps}_{\bar x_i^{\eps}}|,
\end{equation}
which goes to zero as $\eps$ goes to zero by~(\ref{eq:c1b}).

If $t\in[0,t_{k+1}^-]\setminus\mathcal{M}$, then we have that $Z^{(k)}_t\leq\gamma_{x_{k+1}}$ and  
$Z^{(k,\eps)}_{\lambda^\eps(t)}\leq\max_{i>k}\gamma^{\eps}_{x^{\eps}_{i}}$. Hence,  
\begin{equation}
\label{eq:c3}
|Z^{(k)}_t-Z^{(k,\eps)}_{\lambda^\eps(t)}|\leq\gamma_{x_{k+1}}\vee\max_{i>k} 
\gamma^{\eps}_{x^{\eps}_{i}}\leq 
\gamma_{x_{k+1}}\vee\sum_{i>k}\gamma^\eps_{x^{\eps}_{i}}\leq
\sum_{i>k}\gamma_{x_i}+2\eta,
\end{equation}
for all small enough $\eps$, by~(\ref{eq:c1a}).

Now, we solve the first problem by considering two cases:

1) If $s \in [t_i^-, t_i]$ for some $1\leq i\leq k$, then the slope of $\lambda^\eps$ is given by
\begin{equation}
\label{eq:c3a}
\frac{t_i^{\eps}-t_i^{\eps -}}{t_i-t_i^-} = \frac{\gamma^\eps_{\bar x_i^{\eps}} T_i}{ \gamma_{\bar x_i} T_i} 
= \frac{\gamma^\eps_{\bar x_i^{\eps}}}{ \gamma_{\bar x_i}} \rightarrow 1
\end{equation}
as $\eps\to0$, by~(\ref{eq:c1b}).

2) If $s \in [t_i, t_{i+1}^-]$  for some $0\leq i\leq k$, where $t_0:=0$, then it suffices to prove that 
\begin{eqnarray}
\label{eq:c4}
&t_i^{\eps} \rightarrow t_i&\\
\label{eq:c5}
&t_i^{\eps -} \rightarrow t_i^-&
\end{eqnarray}
as $\eps\to0$ in probability.

In all cases, the absolute value of the difference of right and left hand sides is bounded above by
\begin{equation}
\label{eq:c6}
\sum_{i=1}^m|\gamma^\eps_{x_i^{\eps}}-\gamma_{x_i}|T_i+\sum_{i>m}|\gamma^\eps_{x_i^{\eps}}-\gamma_{x_i}|T_i.
\end{equation}

The first term vanishes almost surely as $\eps\to0$ by~(\ref{eq:c1}) and the Law of Large Numbers, 
and, given $\eta>0$, the expected value of the second term is bounded above by
\begin{equation} 
\label{eq:c7}
\sum_{i>m}\gamma^\eps_{x_i^{\eps}}+\sum_{i>m}\gamma_{x_i}\leq2\sum_{i>m}\gamma_{x_i}+\eta,
\end{equation}
for all small enough $\eps$, where use is made of~(\ref{eq:c1a}) in the latter inequality, and~(\ref{eq:c4},\ref{eq:c5}) follow since 
$m\to\infty$ as $\eps\to0$ and $\eta$ is arbitrary.

To conclude, given $0<\zeta<1,\d>0$, choose $T>-\log(\zeta/2)$, and then $R$ satisfying~(\ref{eq:d1}), and then $\d'$  satisfying~(\ref{eq:d2})
with $\eta=\zeta/4$, and then $k$ such that $\sum_{i>k}\gamma_{x_i}<\zeta/2$. Choosing now $\lambda^\eps$ as in~(\ref{eq:lam}), we conclude
that
\begin{equation}
\label{eq:c8}
\limsup_{\eps\to0}\P(d(Z^{(k,\eps)},Z^{(k)})>\zeta)\leq\d,
\end{equation}
where $d$ is the $J_1$ Skorohod distance on $D$ (see~\cite{EK} Chapter 3, Section 5). Since $Z^{(k,\eps)}=Z^\eps$ and $Z^{(k)}=Z$ in distribution 
for all fixed $k$ and $\eps$ small enough, the result follows.
$\square$

Let us now explain how Lemma~\ref{lemma}
will be used in the sequel. Our aim is to apply it to a case where
$S^\eps$ and $S$ are random objects, in fact subordinators, with parameters that are themselves random, which we call {\em environment}. 
Both $S^\eps$ and $S$, as well as their respective environments, will be independent of $\{T_i\}$, and the convergence $S^\eps\rightarrow S$ will hold
only in distribution: either
1) the joint distribution of the environment and the subordinators, 
or 
2) the distributions of subordinators given the environment, for almost every realization of the environment. 
In both cases, we may use the Skorohod representation theorem (see e.g.~\cite{WW} Theorem 3.2.2). In case 1) we will first explicitly choose a
convenient version of the environment, for which the
distribution of the subordinator, given the environment, converges for almost every realization of the environment; 
with the modified environment, we are 
effectively in case 2. We can then, by Skorohod representation, in both cases, for each choice of the environment, choose versions of the 
subordinators that converge almost surely, and then we are in the setting of Lemma~\ref{lemma}. It is clear that the conclusion of the lemma
holds for the original subordinator, where the distribution referred to in the lemma is the joint distribution of $\{T_i\}$ and the subordinators 
given the original environment in case 2, and the modified environment in case 1, for almost every realization of that environment in each case. 
In case 1, 
the result of the lemma will then hold for the overall joint distribution of $\{T_i\}$, the subordinators given the environment, 
and the environment.

Establishing the convergence in distribution of the subordinators is done by verifying the convergence of the respective
Laplace exponents.

\section{Application to trap models on the complete graph and $K$ processes}
\label{sec:tck}
\setcounter{equation}{0}

We will apply the lemma above to show scaling limit results for trap models in the complete graph and for $K$ processes. 
We introduce these two processes next. 

We first consider the trap model on the complete graph 
\begin{equation}
\label{eq:t0}
K_n=\{\{1,\ldots,n\},\,\{(x,y),\,x,y=1,\ldots,n\}\}
\end{equation}
with $n$ vertices (differently from the usual definition, here we include self loops, for convenience -- 
this should not matter in the convergence results below): 
$Y_n=(Y_n(t))_{t\geq 0}$, which is a continuous time Markov chain 
with jump rate at site $x$ given by 
\begin{equation}
\label{eq:t1}
\tau_x^{-(1-a)},
\end{equation}
and transition probability from site $x$ to site $y$ given by
\begin{equation}
\label{eq:t2}
\frac{\tau_y^a}{\sum_{z=1}^n\tau_z^a},
\end{equation}
where $a\in[0,1]$ is a parameter, and 
\begin{equation}
\label{eq:tau}
\tau:=\{\tau_x,\, x=1,2,\ldots\} 
\end{equation}
is an independent family of positive random variables with common distribution
in the domain of attraction of a stable law of degree $0<\alpha<1$, that is,
\begin{equation}
\label{eq:tau1}
\mathbb{P}(\tau_1>t)=\frac{L(t)}{t^{\alpha}},\mbox{     }t>0,
\end{equation}
where $L$ is a slowly varying function at infinity.

We call $Y_n$ an {\em asymmetric} or {\em weighted trap model on the complete graph} with asymmetry parameter $a$, mean jump time parameters
$\{\tau_x^{1-a},\,x=1,\ldots,n\}$ and weights 
$\{\tau_x^{a},\,x=1,\ldots,n\}$. The latter set of parameters may indeed be seen as unnormalized weights of the transition probabilities of $Y_n$.
Notice that the 
$a=0$ (symmetric) case corresponds to uniform weights.

We will consider the following construction of $Y_n$.
Let 
\begin{equation}
\label{eq:caln}
\mathcal{N}=\{N^{(x)}:=(N_r^{(x)})_{r\geq 0},\,x\in\mathbb{N}^*\} 
\end{equation}
be a family of independent Poisson counting processes such that the rate of $N^{(x)}$ is $\tau_x^a$.
Let $\sigma_j^{(x)}$ the $j$-th event time of $N^{(x)}$, $j\geq1$.                                
Let also 
\begin{equation}
\label{eq:calt}
\mathcal{T}=\{T_i^{(x)},\,x\in\mathbb{N}^*\} 
\end{equation}
be independent mean 1 exponential random variables, independent of $\mathcal{N}$ and $\tau$, and define for $r\geq0$
\begin{equation}%{align}
\label{eq:t2a}
S_n(r)=\sum_{x=1}^{n}\tau_x^{1-a}\,N_r^{(x)},\quad\Gamma_n(r)=\sum_{x=1}^{n}\tau_x^{1-a}\sum_{i=1}^{N_r^{(x)}}T_i^{(x)}.
\end{equation}%{align}
Then
\begin{equation}
\label{eq:t3}
 Y_n(t)=x,\mbox{ if }\Gamma_n(\sigma_j^{(x)}-)\leq t<\Gamma_n(\sigma_j^{(x)})\,\mbox{ for some }x,j\geq1.
\end{equation}
is a construction of $Y_n$ as above described, with initial state distributed on $\{1,\ldots,n\}$ in such a way that
site $x$ has probability weight proportional to $\tau_x^a$, $x\in\{1,\ldots,n\}$.

\begin{rmk}
 \label{rmk:id}
Regarding the latter point, notice that the initial state of $Y_n$ is the one whose Poisson mark is the earliest, so it corresponds
to the minimum of $n$ independent exponential random interarrival times with rates $\tau_x^a$, $x\in\{1,\ldots,n\}$, and 
it is well known that the probability that the minimum of $n$ independent exponential random variables is a given such random 
variable is proportional to its rate.
\end{rmk}

Below we will be interested in
\begin{equation}
\label{eq:zn}
 Z_n(t)=\tau^{1-a}_{Y_n(t)}.
\end{equation}

This is the representation for the process aluded to at the introduction above. It has been considered in~\cite{BF}, where the symmetric ($a=0$) case
was studied, and a (single time) scaling limit result was derived for it, first taking the volume, and then the
time, to infinity (see Proposition 2.10 in that reference) -- this is an aging regime not considered in this paper, but rather in~\cite{G}.

\begin{rmk}
 \label{rmk:rezn}
$Z_n$ and $Y_n$ may be seen as processes in random environment, where $\tau$ is the set of random parameters acting as environment. 
Indeed, given $\tau$, both are Markovian (this should be clear for $Y_n$, but a moment's thought reveals that it is true for $Z_n$ 
as well, even when there are same values for $\tau_i$'s with distinct $i$'s). Notice also that $\tau$ is an environment for
$S_n$ as well, which for each $n\geq1$ is a subordinator for every fixed such environment (recall the discussion at the end of Section~\ref{lemma}).
This aspect, which is characteristic of the complete graph, makes our approach particularly suitable, since by an application of (the continuity) 
Lemma~\ref{lemma}, we are left with establishing convergence of subordinators (in the Skorohod topology), which reduces to showing convergence of
Laplace exponents (in the topology of real numbers), which is relatively simple, as we will see below.
\end{rmk}

\begin{rmk}
 \label{rmk:gcon1}
Given $S_n$, $Z_n$ may be identified in distribution to $\Xi(S_n,\{T_i,\,i\in \N^*\})$, with $\Xi$ introduced in Remark~\ref{rmk:int}.
\end{rmk}

We now turn to $K$ processes, which is a Markov process in continuous time on $\bar{\mathbb{N}}^*=\{1,2,\ldots,\infty\}$ constructed in a similar way
as
$Y_n$ was above,
as follows. Let $\g=\{\g_x,\,x\in[0,\infty)\}$ be the increments of an $\a$-stable subordinator in $[0,\infty)$ given by a Poisson process $\cp$ in
$(0,\infty)\times(0,\infty)$ with intensity measure
\begin{equation}
\label{eq:pp} 
\a x^{-1-\a}\,dx\,dy.
\end{equation}

It is well known that the {\em nonzero} set $\{x\in[0,\infty):\,\g_x>0\}$ is countable, so in particular the sums over $[0,1]$ 
below have a countable number of nonzero terms only, and thus make the usual sense, almost surely.

Let 
\begin{equation}
\label{eq:hcaln}
\hat{\mathcal{N}}=\{\hat N^{(x)}:=(\hat N_r^{(x)})_{r\geq 0},\,x\in[0,1]\} 
\end{equation}
be a family of independent Poisson counting processes such that the rate of $\hat N^{(x)}$ is $\g_x^a$,
where $\hat N^{(x)}\equiv0$ whenever $\g_x=0$.
Let $\hat\sigma_j^{(x)}$ the $j$-th event time of $\hat N^{(x)}$, $j\geq1$.
Let also
\begin{equation}
\label{eq:hcalt}
\hat{\mathcal{T}}=\{\hat T_i^{(x)},\,x\in[0,1]\} 
\end{equation}
be a family iid mean 1 exponential random variables independent of $\hat{\mathcal{N}}$.

Define for $r\geq0$
\begin{equation}%{align}
\label{eq:t2ab}
S(r)=\sum_{x\in[0,1]}\g_x^{1-a}\hat N_r^{(x)},\quad\Gamma(r)=\sum_{x\in[0,1]}\g_x^{1-a}\sum_{i=1}^{\hat N_r^{(x)}}\hat T_i^{(x)},
\end{equation}%{align}
and then  make
\begin{equation}
\label{eq:t4}
 Y_t=\begin{cases}
x,&\mbox{ if }\Gamma(\hat\sigma_j^{(x)}-)\leq t<\Gamma(\hat\sigma_j^{(x)})\,\mbox{ for some }x,j\geq1,\\
\infty,&\mbox{ otherwise.}
        \end{cases}
\end{equation}

\begin{rmk}
 \label{rmk:k}
It can be verified that when $a>\a$, then $Y$ is a jump process, and so there is almost surely no $t$ for which $Y(t)=\infty$
(since in this case
$\cup_{j,x}[\Gamma(\hat\sigma_j^{(x)}-),\Gamma(\hat\sigma_j^{(x)}))=[0,\infty)).$
%is a semi-infinite interval. 
And in the case where $a\leq\a$, there almost surely exist $t$'s for which $Y(t)=\infty$.
(One way to check these claims is by verifying that when $a>\a$, $\{\hat\sigma_j^{(x)};\,j\geq1, x\in[0,1]\}$ 
is a discrete subset of $[0,\infty)$ almost surely, and when $a\leq\a$, it is almost surely dense in
$[0,\infty)$, and these in turn follow from the fact that $\sum_{x\in[0,1]}\g_x^a$ is almost surely finite in the 
former case, and infinite in the latter one.)
\end{rmk}

Let
\begin{equation}
\label{eq:zt}
 Z_t=\gamma^{1-a}_{Y_t},
\end{equation}
where $\gamma_\infty$ should be interpreted as $0$.
\begin{rmk}
 \label{rmk:rez}
$Z$ and $Y$ may be seen as processes in random environment, where $\gamma$ (more specifically, $\gamma|_{[0,1]}=\{\g_x,\,x\in[0,1]\}$) 
is the environment. Indeed, given $\g$, both are Markovian. $\gamma|_{[0,1]}$ is also an environment for
$S$, which is a subordinator for every fixed such environment (recall the discussion at the end of Section~\ref{lemma}).
\end{rmk}

\begin{rmk}
 \label{rmk:gcon2}
Given $S$, $Z$ may be identified in distribution to $\Xi(S,\{T_i,\,i\in \N^*\})$, with $\Xi$ introduced in Remark~\ref{rmk:int}.
\end{rmk}

\begin{rmk}
 \label{rmk:rep}
In~\cite{FM1} and other references the representations used for the trap model and $K$ process are the ones given here by $Y_n(t)$ and $Y(t)$,
$t\geq0$,
respectively (see~(\ref{eq:t3}) and~(\ref{eq:t4}) above). The alternative representation $Z_n(t)$ and $Z(t)$, $t\geq0$, we adopt here
(see~(\ref{eq:zn}) and~(\ref{eq:zt}) above) has
the advantage of leading to a unifying approach for taking the scaling limits of those processes, as explained in the introduction and will be done
in detail in Subsections~\ref{ssec:erg}
and~\ref{ssec:agez} below.
\end{rmk}

In the next subsection, we will consider a particular scaling regime for $Z_n$ and establish a scaling limit result under which $Z_n$ converges to the
$K$ process. Then, in the following subsection
we will derive a scaling limit result satisfied by $Z$.
All proofs will rely on Lemma~\ref{lemma} above to get the results from the convergence of the appropriate $S^\eps$ in each case (see statement of
that lemma and its preliminaries above). In order to obtain the latter convergence, since we have subordinators in all cases, it will suffice to
establish convergence of the associated Laplace exponents. The last subsection is devoted to a discussion on aging results (for two-time correlation
functions) satisfied by $Z$ as a consequence of Theorem~\ref{teo:asz} and other results.

\subsection{Scaling limit for $Z_n$ at large times}
\label{ssec:erg}
For $r\geq0$, let
\begin{equation}\label{eq:u}
U(r)=\sum_{x\in[0,r]}\g_x.
\end{equation}

Given a sequence $(c_n)_{n\geq1}$, set
\begin{equation}
\label{eq:sczn}
Z_t^{(n)}=c_n^{1-a} Z_n(t/c_n^{1-a}),\,\,t\geq0.
\end{equation}

Let $P_1$ denote the probability measure induced on $(D,J_1)$ by the joint distribution of $\tau$,
${\cal N}$ and ${\cal T}$ -- given above in respectively~(\ref{eq:tau}),~(\ref{eq:caln}) and~(\ref{eq:calt}).

\begin{theorem}
\label{teo:es}
There exists a deterministic sequence $(c_n)_{n\geq1}$ such that
\begin{equation}
\label{eq:znz}
(Z_t^{(n)})_{t\geq0}\joa(Z_t)_{t\geq0}.
\end{equation}
as $n\to\infty$.
\end{theorem}

The sequence $(c_n)$ will be exhibited explicitly in the proof below (see~\ref{eq:cn}).

 \begin{proof}\mbox{}

By Lemma~\ref{lemma}, and recalling the discussion at the end of Section~\ref{sec:lm}, it is enough to establish the limit
\begin{equation}
\label{eq:e1}
S^{(n)}\joa S,
\end{equation}
where
\begin{equation}
\label{eq:e1a}
S_r^{(n)}:=c_n^{1-a}S_n(c_n^ar)=\sum_{x=1}^n(c_n\tau_x)^{1-a}N_{c_n^ar}^{(x)},
\end{equation}
since, given $S^{(n)}$, $Z^{(n)}$ is identically distributed with $\Xi(S^{(n)},\{T_i,\,i\in \N^*\})$ -- see Remark~\ref{rmk:int}. 

In order to establish~(\ref{eq:e1}), we will make a precise choice of $c_n$ and switch to another
version of $\tau$, which properly rescaled converges strongly, rather than weakly.
We follow \cite{FIN}, Section 3. Let
\begin{eqnarray}
  \label{eq:cn}
  &c_n=\left(\inf\{t\geq0:\P(\tau_1>t)\leq n^{-1}\}\right)^{-1},&\\
  \label{eq:taun}
  &\tau_x^{(n)}:=c_n^{-1}\,g_n\!\left(U(x)-U(x-1/n)\right),\,x\in(0,1]\cap\frac1n\Z&\\
  \label{eq:gn}
  &g_n(y)=c_n\,G^{-1}(n^{1/\a}y),\,y\geq0,
\end{eqnarray}
where $G^{-1}$ is the inverse of the function $G$ defined by the following condition.
\begin{equation}
\label{eq:G}
\P(U(1)>G(x))=\P(\tau_1>x),\,x\geq0
\end{equation}

We then have that $\tau^{(n)}:=\{\tau_x^{(n)},\,x\geq1\}$ is equally distributed with $\tau$ for every $n\geq1$.

For $x\in(0,1]\cap\frac1n\Z$, let now
\begin{equation}
\label{eq:gan}
\g^{(n)}_x=c_n\tau_x^{(n)},
\end{equation}
and define
\begin{equation}
\label{eq:tsn}
\tilde S_r^{(n)}:=\sum_{x=1}^n(\g^{(n)}_{x/n})^{1-a}\tilde N_{r}^{(n,x)},
\end{equation}
where, given $\g$,
\begin{equation}
\label{eq:tcaln}
\tilde{\mathcal{N}}^{(n)}=\{\tilde N^{(n,x)}:=(\tilde N_r^{(n,x)})_{r\geq 0},\,x\in\N^\ast\} 
\end{equation}
is a family of independent Poisson counting processes such that the rate of $N^{(n,x)}$ is $(\g^{(n)}_x)^a$.

One now readily checks, using the identity in distribution of $\tau^{(n)}$ and $\tau$ for every $n\geq1$, together with the above definitions, 
that $\tilde S^{(n)}:=(\tilde S_r^{(n)})_{r\geq0}$ has the same distribution (induced by $(\g,\tilde{\mathcal{N}}^{(n)})$) as $S^{(n)}$ 
under $P'$ for every $n\geq1$. So it is enough to show that
\begin{equation}
\label{eq:e1c}
\tilde S^{(n)}\jta S,
\end{equation}
where $P_2$ is the probability measure induced on $(D,J_1)$ by the joint distribution of $\g$ and $\tilde{\mathcal{N}}^{(n)}$.

Now since, given $\g$, $\tilde S^{(n)}$ is a subordinator for each $n\geq1$, it is enough to show the
convergence of the Laplace exponents of $\tilde S^{(n)}$, $n\geq1$, as $n\to\infty$, for almost every realization of $\g$,
to the Laplece exponent of $S$ given $\g$, which is itself a subordinator. (See Corollary 3.6 page 374 in \cite{JS}.)

A straightforward computation yields
\begin{equation}
\label{eq:e3}
\tilde\varphi_n(\l):=\sum_{x\in(0,1]\cap\frac1n\Z}(\g_x^{(n)})^a(1-e^{-\l (\g_x^{(n)})^{1-a}})
\end{equation}
as the Laplace exponent of $\tilde S^{(n)}$, $n\geq1$.

Now let
\begin{equation}
\label{eq:taudelta}
\t_{\delta}=\{x\in[0,1]:\, \g(x)>\delta\}=\{x_1<\ldots<x_K\}, 
\end{equation}
and 
\begin{equation}
\label{eq:taudeltan}
\t^{(n)}_{\delta}=\left\{x^{(n)}_1=\frac1n\lceil nx_1\rceil<\ldots<x^{(n)}_K=\frac1n\lceil nx_K\rceil\right\}, 
\end{equation}
where the strict inequalities in~(\ref{eq:taudeltan}) hold provided $n$ is large enough (for each fixed $\d$).

Lemma 3.1 in~\cite{FIN} implies that for every $\d>0$
\begin{equation}
\label{eq:e4}
\sum_{x\in\t^{(n)}_{\delta}}(\g_x^{(n)})^a(1-e^{-\l (\g_x^{(n)})^{1-a}})\to\sum_{x\in\t_{\delta}}\g_x^a(1-e^{-\l \g_x^{1-a}})
\end{equation}
almost surely as $n\to\infty$. One also readily checks that
\begin{equation}
\label{eq:e5}
\sum_{x\in(0,1]\cap\frac1n\Z\setminus\t^{(n)}_{\delta}}(\g_x^{(n)})^a(1-e^{-\l (\g_x^{(n)})^{1-a}})
\leq \l\sum_{x\in(0,1]\cap\frac1n\Z\setminus\t^{(n)}_{\delta}}\g_x^{(n)}.
\end{equation}
Since, as argued in paragraphs of (3.25-3.28) in~\cite{FIN}, we have that the $\lim_{\delta\to0}\limsup_{n\to\infty}$ of the sum in the right hand
side of~(\ref{eq:e5}) vanishes almost surely, we may conclude that
\begin{equation}
\label{eq:e6}
\tilde\varphi_n(\l)\to\varphi(\l):=\sum_{x\in[0,1]}\g_x^a(1-e^{-\l \g_x^{1-a}}),\,\,\l\geq0,
\end{equation}
almost surely. This convergence holds in principle for each $\l\geq0$, but it may be argued to hold simultaneously for every $\l\geq0$
from the monotonicity of $\tilde\varphi_n$ for every $n\geq1$, and the continuity of $\varphi$.  The right hand
side of~(\ref{eq:e6}) is the Laplace exponent of $S$ given $\g$, so the proof is complete.
\end{proof}

\subsection{Scaling limit of $Z$ at small times}
\label{ssec:agez}

In this subsection, we assume $0\leq a<\a$.
Let
\begin{equation}
\label{eq:scz}
Z_t^{(\eps)}=\eps^{-1} Z_{\eps t}.
\end{equation}

Before stating a convergence result for $Z^{(\eps)}$, let us describe the limit process.
Let $(\hat S_t)_{t\geq0}$ be an $\hat\a$-stable subordinator, where 
\begin{equation}
\label{eq:ha}
\hat\a=\frac{\alpha -a}{1-a}, 
\end{equation}
and whose Laplace exponent is given by 
$\hat\varphi(\l)=\hat c\l^{\hat\a}$, where $\hat c$ is a constant to be determined below. %~(\ref{eq:hc})). 

We may then write $\hat S$ as a partial sum of its increments as follows.
\begin{equation}
\label{eq:hu}
\hat S_r=\sum_{x\in[0,r]}\hat\g_x,
\end{equation}
where $\{\hat\g_x,\,x\in\N^*\}$ are the increments of $\hat S$.

Let now 
\begin{equation}
\label{eq:hv}
\hat\G_r=\sum_{x\in[0,r]}\hat\g_xT_x,
\end{equation}
where 
\begin{equation}
\label{eq:ctp}
{\cal T}':=\{T_x,\,x\in[0,\infty)\} 
\end{equation}
is an iid family of mean 1 exponential random variables, independent of $\hat S$. 

\begin{rmk}
 \label{rmk:v}
One may readily check that $\hat\G$ is also an $\hat\a$-stable subordinator (under the joint distribution of $\hat S$ and
$\{T_x,\,x\in[0,\infty)\}$).
\end{rmk}

Now define
\begin{equation}
\label{eq:hz}
\hat Z_t=\begin{cases}
              \hat\gamma_{x},&\mbox{  if }t\in[\hat\Gamma_{x-},\hat\Gamma_{x})\mbox{ for some }x\in[0,\infty)\\
\mbox{}\,\,\,0,&\mbox{  for all other }t\geq0,\mbox{  if any}.
             \end{cases}
\end{equation}

\begin{rmk}
 \label{rmk:rehz}
$\hat Z$ may be seen as a process in random environment, where $\hat S$ is the environment. Indeed, given $\hat S$, $\hat Z$ is Markovian.
And the distribution of  $\hat Z$ (integrated over the environment) makes it a self similar process of index $1$, that is,
$(\hat Z_t)_{t\geq0}=(c^{-1}\hat Z_{ct})_{t\geq0}$ in distribution for every constant $c>0$. This latter property explains the aging behavior of
$Z$ in its small time scaling regime, as established below.
\end{rmk}

\begin{rmk}
 \label{rmk:gcon3}
Given $\hat S$, $\hat Z$ may be identified in distribution to $\Xi(\hat S,\{T_i,\,i\in \N^*\})$, with $\Xi$ introduced in Remark~\ref{rmk:int}.
\end{rmk}

Before we state this subsection's result, let, for $\g$ fixed, $P_3=P_3^\g$ denote the the probability measure induced on 
$(D,J_1)$ by the joint distribution of $\hat{\cal N}$ and $\hat{\cal T}$ -- given above in respectively~(\ref{eq:hcaln})
and~(\ref{eq:hcalt}).

\begin{theorem}
\label{teo:asz}
If $0\leq a<\a$ then for almost every $\g$
\begin{equation}
\label{eq:zu}
(Z_t^{(\eps)})_{t\geq0}\jtha(\hat Z_t)_{t\geq0}.
\end{equation}
as $\eps\to0$.
\end{theorem}

\begin{rmk}
 \label{rmk7}
Perhaps more precisely, Theorem~\ref{teo:asz} states that for almost every $\g$, the distribution of 
$(Z_t^{(\eps)})$ under $P_3$ converges to that of $(\hat Z_t)$ under $P_4$, the probability measure induced on 
$(D,J_1)$ by the joint distribution of $\g$ and ${\cal T}'$.% via the construction leading to~(\ref{eq:hz}).
\end{rmk}

\begin{cor}
\label{cor}
If $0\leq a<\a$ then 
\begin{equation}
\label{eq:zu1}
(Z_t^{(\eps)})_{t\geq0}\jfia(\hat Z_t)_{t\geq0}
\end{equation}
as $\eps\to0$, where $P_5$ denotes the probability measure induced on 
$(D,J_1)$ by the joint distribution of $\g$,~$\hat{\cal N}$ and $\hat{\cal T}$.
\end{cor}

\begin{rmk}
 \label{rmk7+}
The above corollary follows immediately from the preceding theorem, since $P_5$ is obtained by integrating $P_3$ over
the distribution of $\g$. Below we will nevertheless give a direct (sketchy) argument for the corollary, much simpler 
than the one for the theorem next.
\end{rmk}

{\noindent\bf Proof of Theorem~\ref{teo:asz}}

Let
\begin{equation}
\label{eq:se}
\hat S_r^{(\eps)}=\eps^{-1}\sum_{x\in[0,1]}\gamma_x^{1-a}\hat N^x_{\eps^{\hat\a}r},\,r\geq0
\end{equation}
where $\hat\a$ was introduced in~(\ref{eq:ha}) above.
Then, given $\g$ and $\eps>0$, $(\hat S_t^{(\eps)},\,t\geq0)$ is a subordinator, and its Laplace exponent equals
\begin{equation}
\label{eq:lese}
\hat\varphi^{(\eps)}(\l)=\eps^{\hat\a}\sum_{x\in[0,1]}\gamma_x^a(1-e^{-\l\eps^{-1}\gamma_x^{1-a}}),\,\,\l\geq0.
\end{equation}

By Lemma~\ref{lemma}, and recalling the discussion at the end of Section~\ref{sec:lm}, to get the result, it is enough to establish the limit
\begin{equation}
\label{eq:e1b}
\hat S^{(\eps)}\jfia\hat S
\end{equation}
%in distribution 
as $\eps\to0$ for a.e.~$\g$.
Since we are dealing with subordinators, it suffices to show that for almost every $\g$
\begin{equation}
\label{eq:a1}
\hat\varphi^{(\eps)}(\l)\to\hat c\l^{\hat\a},\,\,\l\geq0,
\end{equation}
as $\eps\to0$, for some positive finite constant $\hat c$.
This is obvious for $\l=0$, so let us fix $\l>0$, and write
\begin{equation}
\label{eq:a2}
\l^{-\hat\a}\hat\varphi^{(\eps)}(\lambda)=R^{-\alpha}\sum_{x\in[0,1]}(R\gamma_x)^a(1-e^{-(R\gamma_x)^{1-a}})
\end{equation}
with $R=(\eps^{-1}\lambda)^{\frac{1}{1-a}}$, and then argue in the sequel that the left hand side converges to a constant as $R\to\infty$
for a.e.~$\gamma$.

We start by considering
\begin{equation}
\label{eq:y}
W:=R^{-\alpha}\sum_{x\in
[0,1]}\sum_{i=1}^{R\delta^{-1}}(R\gamma_x)^a(1-e^{-(R\gamma_x)^{1-a}})\mathbb{I}_{\{\gamma_x\in[\frac{\delta}{R}(i-1),\frac{\delta}{R}i]\}}.
\end{equation}
Since the difference between $W$ and the left hand side of~(\ref{eq:a2}) is bounded above by
\begin{equation}
R^{-(\alpha-a)}\sum_{x\in[0,1]}\gamma_x^a\,\mathbb{I}_{\{\gamma_x>1\}}, 
\end{equation}
which vanishes as $R\to\infty$ for  a.e.~$\gamma$, it is enough to establish the
convergence result for $W$. We estimate it as follows.
\begin{eqnarray}
\label{E:prim}
W-X_1&\leq& R^{-\alpha}\sum_{i=2}^{R\delta^{-1}}X_i^+:=R^{-\alpha}\sum_{i=2}^{R\delta^{-1}}(\delta i)^a(1-e^{-(\delta i)^{1-a}})M_i\\
\label{E:sec}
W&\geq& R^{-\alpha}\sum_{i=2}^{R\delta^{-1}}X_i^-:=R^{-\alpha}\sum_{i=2}^{R\delta^{-1}}(\delta (i-1))^a(1-e^{-(\delta (i-1))^{1-a}})M_i,
\end{eqnarray}
where $X_1=R^{-\alpha}\sum_{x\in[0,1]}(R\gamma_x)^a(1-e^{-(R\gamma_x)^{1-a}})
\mathbb{I}_{\{\gamma_x\in[0,\frac{\delta}{R}]\}}$ and $M_i$ is the number of points of $\cp$ in the region
$[0,1]\times[\frac{\delta}{R}(i-1),\frac{\delta}{R}i]$ (recall paragraph of~(\ref{eq:pp}) above). 

$X_1$ can be bounded above by $R^{-\alpha}\sum_{x\in[0,1]}(R\gamma_x)\,
\mathbb{I}_{\{\gamma_x\in[0,\frac{\delta}{R}]\}}$, and this has the same distribution as 
$R^{-\alpha}\sum_{x\in[0,R^\a]}\gamma_x\,\mathbb{I}_{\{\gamma_x\in[0,\delta]\}}$ for every $R>0$, 
by the scale invariance of $\g$. We can use standard large deviation estimates 
for the latter expression to conclude that $X_1$ can be ignored in the limits as $R\to\infty$ and then $\d\to0$ (here we may use the existence of a
positive 
exponential moment for $\sum_{x\in[0,1]}\gamma_x\,\mathbb{I}_{\{\gamma_x\in[0,\delta]\}}$ for any $\d$, a result that follows as an application of
Campbell Theorem -- see~\cite{K}). 
We concentrate on the right hand sides of~(\ref{E:prim}, \ref{E:sec}).

We start with~(\ref{E:prim}). By the exponential Markov inequality, we get, for given $\theta,\xi>0$,
\begin{equation}
\label{eq:a3}
\mathbb{P}\left(R^{-\alpha}\sum_{i=2}^{R\delta^{-1}}X_i^+\geq R^{-\alpha}\sum_{i=2}^{R\delta^{-1}}\mathbb{E}X^+_i+\xi\right)\leq
\frac{A}{B}%\nonumber
\end{equation}
where $A=\mathbb{E}e^{\theta\sum_{i=1}^{R\delta^{-1}}X_i^+}$ and $B=e^{\theta\sum_{i=1}^{R\delta^{-1}}\mathbb{E}X^+_i+ R^{\alpha}\xi }$. 

Since $M_i$, $i\geq2$, are independent Poisson random variables, we obtain
\begin{equation}
\label{eq:a4}
\frac{A}{B}=e^{-R^{\alpha}\xi\theta+\sum_{i=2}^{R\delta^{-1}}(e^{c_i\theta}-1-c_i\theta)\,\mathbb{E}M_i},
\end{equation}
where $c_i=(\delta i)^a(1-e^{-(\delta i)^{1-a}})$.

We choose $\theta=R^{-b}$ with $a<b<\alpha<2b$. Then, using the estimate
\begin{equation}
\label{eq:a5}
\mathbb EM_i=\int_{\frac{\delta}{R}(i-1)}^{\frac{\delta}{R}i}\frac{\alpha}{x^{1+\alpha}}dx\leq\frac{R^{\alpha}}{\delta^{\alpha}(i-1)^{1+\alpha}},
\end{equation}
we find that the sum in the exponent in~(\ref{eq:a4}) is bounded above by
\begin{equation}
\label{eq:a6}
\sum_{i=2}^{R\delta^{-1}}\frac{R^{\alpha}}{\delta^{\alpha}(i-1)^{1+\alpha}}(c_iR^{-b})^2
\leq2\frac{R^{\alpha-2b}}{\delta^{\alpha-2a}}\sum_{i=1}^{R\delta^{-1}}i^{-(1+\alpha-2a)}%(1-e^{-(\delta i)^{1-a}})^2}.
\end{equation}
Since the sum on the right of~(\ref{eq:a6}) is bounded by constant times $R^{2a-\a}\vee\log R$, and using the above estimates, we 
find that the exponent in~(\ref{eq:a4}) is bounded above by
\begin{equation}
\label{eq:a7}
-R^{\alpha-b}\xi+\mbox{ const } R^{-c'},
\end{equation}
for some constant $c'>0$. We can then apply Borel-Cantelli and conclude that for a.e.~$\g$, given $\xi>0$
\begin{equation}
\label{eq:a8}
R^{-\alpha}\sum_{i=2}^{R\delta^{-1}}X_i^+\leq R^{-\alpha}\sum_{i=2}^{R\delta^{-1}}\mathbb{E}X^+_i +\xi
\end{equation}
for all large enough  $R$.

Conversely, we can conclude that given $\xi>0$, for a.e.~$\g$ and all $R$ large enough
\begin{equation}
\label{eq:a9}
R^{-\alpha}\sum_{i=2}^{R\delta^{-1}}X_i^-\geq R^{-\alpha}\sum_{i=2}^{R\delta^{-1}}\mathbb{E}X^-_i-\xi.
\end{equation}
(\ref{eq:a8}) and (\ref{eq:a9}) then imply that
\begin{equation}
\label{eq:a10}
\liminf_{R\to\infty} R^{-\alpha}\sum_{i=1}^{R\delta^{-1}}\mathbb{E}X^-_i\leq\liminf_{R\to\infty} W
\leq\limsup_{R\to\infty} W\leq \limsup_{R\to\infty}  R^{-\alpha}\sum_{i=1}^{R\delta^{-1}}\mathbb{E}X^+_i.
\end{equation}

To conclude, it is enough to verify that
\begin{equation}
\label{limite}
\liminf_{\delta\rightarrow 0}\liminf_{R\to\infty} R^{-\alpha}\sum_{i=1}^{R\delta^{-1}}\mathbb{E}X^-_i
=\limsup_{\delta\rightarrow 0}\limsup_{R\to\infty}R^{-\alpha}\sum_{i=1}^{R\delta^{-1}}\mathbb{E}X^+_i
\end{equation}
is a (positive finite) constant $\hat c$.

We begin with the following estimate.
\begin{equation}
\mathbb{E}X_i^+=(\delta i)^a(1-e^{-(\delta i)^{1-a}})\int_{\frac{\delta}{R}(i-1)}^{\frac{\delta}{R}i}\frac{\a}{x^{1+\alpha}}dx
\leq (\delta i)^a(1-e^{-(\delta i)^{1-a}})\frac{\delta}{R}\frac{\a}{(\frac{\delta}{R}(i-1))^{1+\alpha}}
\end{equation}
Summing up:
\begin{eqnarray}\nonumber
R^{-\alpha}\sum_{i=2}^{R\delta^{-1}}\mathbb{E}X_i^+
&\leq& R^{-\alpha}\sum_{i=2}^{R\delta^{-1}}(\delta i)^a(1-e^{-(\delta i)^{1-a}})\frac{\delta}{R}\frac{\a}{(\frac{\delta}{R}(i-1))^{1+\alpha}}\\
&=&\a\sum_{i=2}^{R\delta^{-1}}\delta^{a-\alpha}\frac{1-e^{-(\delta i)^{1-a}}}{i^{1+\alpha-a}}\left(\frac i{i-1}\right)^{1+\a}
\end{eqnarray}
Now as $R\to\infty$, the latter sum converges to a series, which is readily seen to be an approximation to an integral. We find that
\begin{equation}
\label{eq:a11}
\limsup_{\delta\rightarrow 0}\limsup_{R\to\infty} R^{-\alpha}\sum_{i=2}^{R\delta^{-1}}\mathbb{E}X_i^+\leq
\a\int_0^{\infty}\frac{1-e^{-x^{1-a}}}{x^{1+\alpha-a}}dx.
\end{equation}
We similarly find the latter expression as a lower bound for 
$$\liminf_{\delta\rightarrow 0}\liminf_{R\to\infty} R^{-\alpha}\sum_{i=2}^{R\delta^{-1}}\mathbb{E}X_i^-$$
and~(\ref{limite}) follows, with the right hand side of~(\ref{eq:a11}) as the constant $\hat c$. $\square$

\medskip

{\noindent\bf  (Direct) Proof of Corollary~\ref{cor}} (sketchy)

Under $P_5$ we may use a different, more suitable version of $\g$.
In view of the right hand side of~(\ref{eq:a2}), we replace $\g_x$ by $R^{-1}\gamma_{R^\a x}$, $x\in[0,1]$, with $R$ as in the above proof.
The Laplace exponent of the corresponding version of $\hat S^{(\eps)}$ is then readily seen to equal
\begin{equation}
\label{eq:a2p}
\l^{\hat\a}\,R^{-\alpha}\!\!\!\!\sum_{x\in[0,R^{\alpha}]}\gamma_x^a(1-e^{-\gamma_x^{1-a}}).
\end{equation}
Since $\sum_{x\in[0,1]}\gamma_x^a(1-e^{-\gamma_x^{1-a}})$ is integrable, with mean $\hat c$, as can be checked by an 
application of Campbell Theorem, the Law of Large Numbers yields the almost sure convergence of~(\ref{eq:a2p}) to
$\hat c\l^{\hat\a}$ (simultaneously for all $\l\geq0$, once one uses monotonicity and continuity of the functions involved,
as previously argued -- see the end of the proof of Theorem~\ref{teo:es} above). Since that is the Laplace exponent of 
$(\hat Z_t)$, we conclude that the version of $(Z_t^{(\eps)})$ with $\g$ replaced by
$\g^{R}:=\{R^{-1}\gamma_{R^\a x},\,x\in[0,1]\}$ converges in $P_3^{\g^{R}}$-distribution to $(\hat Z_t)$ for almost 
every $\g$. Upon integrating over the distribution of $\g$, we get the convergence in $P_5$-distribution. $\square$

\begin{rmk}
 \label{rmk4}
A few words about the cases where $\a\leq a\leq1$. When $a>\a$, we have that $Z$ is a jump process
in $\N^*$ (see Remark~\ref{rmk:k}) with $Z(0)=\g_x$ with probability proportional to $\g_x^a$, $x\in[0,1]$. It follows 
then that  $Z^{(\eps)}\to\infty$ identically almost surely as $\eps\to0$. 

The case $a=\a$ demands more delicate analysis. We have that $\hat\varphi^{(\eps)}(\l)$ (see~\ref{eq:lese}), 
when scaled with a factor of $|\log\eps|^{-1}$ (instead of $\eps^{\hat\a}=1$ in this case), 
converges to a number $r$ independent of $\l>0$ as $n\to\infty$ in probability, and this is the Laplace exponent of a
subordinator which equals $0$ for an exponentially distributed amount of time of rate $r$, and then jumps to $\infty$, 
where it stays. One may then argue from this that $Z^{(\eps)}\to\infty$ identically as $\eps\to0$ in 
probability.
\end{rmk}

\subsection{Aging in the $K$ process}
\label{ssec:agezz}
Theorem~\ref{teo:asz} may be viewed as an aging result for $Z$, since $\hat Z$ is nontrivial and self similar with
index 1. Corresponding aging results for two-time correlation functions follow.

Below we consider three examples of correlation functions related to aging, and derive scaling limit/aging results for them as a consequence of
Theorem~\ref{teo:asz} (as well as of other results derived above). 
Other correlation functions can be similarly treated. 

\paragraph{Example 1}

We start with the time correlation function introduced in~\cite{B}, which is the one that is usually studied 
in connection with his model. Let
\begin{equation}
\label{eq:ag2a}
\bar\Pi(t,s;\g)=\P(\mbox{no jump of } Z\mbox{ on }[t,t+s]|\g)
\end{equation}
(see Remark~\ref{rmk8} below).

Let $\Phi\in D$, let $\D(\Phi)$ denote the set of discontinuities of $\Phi$, that is,
$\D(\Phi)=\{t\geq0:\Phi(t)\ne \Phi(t-)\}$, and consider $F:D\times(0,\infty)\times(0,\infty)\to\{0,1\}$ such that
\begin{equation}
\label{eq:ag1}
F(\Phi;t,s)=1\{[t,t+s]\cap\D(\Phi)=\emptyset\}.
\end{equation}
Then we have that
\begin{equation}
\label{eq:ag3}
\bar\Pi(\eps t,\eps s;\g)=\esp[F(Z^{(\eps)};t,s)|\g].
\end{equation}
Let also
\begin{equation}
\label{eq:ag2b}
\hat\Pi(t,s)=\P(\mbox{no jump of } \hat Z\mbox{ on }[t,t+s]).
\end{equation}

Since deterministic single times are almost surely continuity points of $\hat Z$, we have that $F(\cdot;t,s)$ is almost surely continuous
under the distribution of $\hat Z$. We thus conclude from
Theorem~\ref{teo:asz} that if $0\leq a<\a$, then for almost every $\g$
\begin{equation}
\label{eq:ag4}
\lim_{\eps\to0}\bar\Pi(\eps t,\eps s;\g)=\esp[F(\hat Z;t,s)]=\hat\Pi(t,s).
\end{equation}

The aging phenomenon, namely $\hat\Pi(\cdot,\cdot)$ being a (nontrivial) function of the ratio of its arguments, then follows from the self similarity
with index 1 (and nontriviality) of $\bar Z$, but in this case there is an explicit expression for $\hat\Pi$, obtained as follows.
One readily checks that the right hand side of~(\ref{eq:ag2b}) equals $\P([t,t+s]\cap\RR(\hat\G)=\emptyset)$, where $\RR(\Phi)$ is the range of
$\Phi\in D$.
Since $\hat\G$ is an $\hat\a$-stable subordinator (see Remark~\ref{rmk:v}), an application of the Dynkin and Lamperti arcsine law theorem for that
probability yields
\begin{equation}
\label{eq:ag5}
\hat\Pi(t,s)=\frac{\sin(\pi\hat\alpha)}{\pi}\int_{s/(t+s)}^1\theta^{-\hat\a}(1-\theta)^{\hat\a-1}\,d\theta.
\end{equation}

The limit in~(\ref{eq:ag4}) was first obtained in~\cite{B} (as the expression in~(\ref{eq:ag5})) for the case where $a=0$. The general case $0\leq
a\leq1$ was first studied in~\cite{G} (see Theorem 3.3 for the case $a < \alpha$, and Theorem 3.4 for the case  $a > \alpha$; the particular
limit~(\ref{eq:ag4}) is (7.5) in that reference).

In case $a\geq\a$, then the discussion in Remark~\ref{rmk4} %and~\ref{rmk5} 
indicates that the limit in~(\ref{eq:ag4}) 
is identically 1, and that aging is thus interrupted.

For the next examples, we restrict $a$ to $[0,\a)$.

\paragraph{Example 2} 

Let
\begin{equation}
\label{eq:ag6}
\bar R(t,s;\gamma)=\P(Z(t)=Z(t+s)|\gamma).
\end{equation}
Then, the difference between $\bar R(\eps t,\eps s;\gamma)=\P(\hat Z^{(\eps)}_t=\hat Z^{(\eps)}_{t+s}|\gamma)$ and
$\bar\Pi(\eps t,\eps s;\gamma)$
is given by
\begin{equation}
\label{eq:ag7}
\P(\hat Z^{(\eps)}_t=\hat Z^{(\eps)}_{t+s};\,\hat Z^{(\eps)}_t\ne \hat Z^{(\eps)}_{t+r}\mbox{ for some }r\in[0,s]|\gamma).
\end{equation}
Let $\hat\p^{(\eps)}$ and $\hat\p$ denote the point processes in $(0,\infty)\times(0,\infty)$ 
associated to $\hat S^{(\eps)}$ and $\hat S$, respectively, i.e., 
\begin{equation}
\label{eq:pps}
\hat\p^{(\eps)}=\left\{\!\left(t,\hat S^{(\eps)}_t-\hat S^{(\eps)}_{t-}\right)\!\!:\,t>0,\,\hat S^{(\eps)}_t-\hat S^{(\eps)}_{t-}>0\right\},\,
\hat\p=\left\{\!\left(t,\hat S_t-\hat S_{t-}\right)\!\!:\,t>0,\,\hat S_t-\hat S_{t-}>0\right\}.
\end{equation}
The convergence in distribution 
$\hat S^{(\eps)}\to\hat S$ argued in the proof of Theorem~\ref{teo:asz} implies that 
\begin{equation}
\label{eq:ag8}
\hat\p^{(\eps)}\to\hat\p
\end{equation}
as $\eps\to0$ in distribution (in the point process sense; for almost every $\gamma$). 

Let also
$\hat\G^{(\eps)}(t)=\eps^{-1}\G(\eps^{\hat\a} t)$, $t\geq0$ (see paragraph of~(\ref{eq:t2ab}) above). 
We have that
\begin{equation}
\label{eq:ag10}
\hat \G^{(\eps)}\to\hat\G
\end{equation}
in distribution for almost every $\gamma$ (see~(\ref{eq:hv}) above). This claim may be argued as follows.
Since $(\hat \G^{(\eps)}_t)$ is a subordinator, an entirely similar reasoning to the one employed in the proof of Theorem~\ref{teo:asz}
may also be employed to establish this result. It also follows from a continuity property of 
$(\hat \G^{(\eps)}_t)$ as a function of $(\hat S^{(\eps)}_t)$ and $\mathcal{T}$ similar to the one established in Lemma~\ref{lemma},
and similarly proven. We leave the details for the interested reader.

For arbitrary $\d,T>0$, consider now the event
\begin{equation}
\label{eq:ag9}
A^{(\eps)}_{\d,T,t,s}=\{\hat Z^{(\eps)}_t>\d,\,\hat Z^{(\eps)}_{t+s}>\d,\,\hat\G^{(\eps)}(T)>t+s\},%\cap A^{(\eps)}_{\d,T},
\end{equation}
and let $B^{(\eps)}_{\d,T}$ be the event that there exist two points in $\hat\p^{(\eps)}\cap\{(0,2T)\times(\d/2,\infty)\}$ with the same 
second coordinate. Now one readily gets from the above convergence results that
\begin{equation}
\label{eq:ag11}
\lim_{\eps\to0}\P(A^{(\eps)}_{\d,T,t,s}|\gamma)=\P(\hat Z_t>\d,\,\hat Z_{t+s}>\d,\,\hat\G(T)>t+s),
\end{equation}
and this can be made arbitrarily close to 1 by choosing $\d$ and $T$ appropriately. We also have that
\begin{equation}
\label{eq:ag12}
\lim_{\eps\to0}\P(B^{(\eps)}_{\d,T}|\gamma)=\P(B_{\d,T}),
\end{equation}
where $B_{\d,T}$ is the event corresponding to $B^{(\eps)}_{\d,T}$ upon replacing $\hat\p^{(\eps)}$ by $\hat\p$.
The latter probability clearly vanishes. Now, since the intersection of the event in the probability in~(\ref{eq:ag7}) 
and $A^{(\eps)}_{\d,T,t,s}$ is contained in $B^{(\eps)}_{\d,T}$, we conclude from the above that
\begin{equation}
\label{eq:ag13}
\lim_{\eps\to0}\bar R(\eps t,\eps s;\gamma)=\lim_{\eps\to0}\bar\Pi(\eps t,\eps s;\gamma)=\hat\Pi(t,s)
\end{equation}
for almost every $\gamma$.

\begin{rmk}
 \label{rmk8}
The {\em aging} correlation functions
\begin{eqnarray}
\label{eq:ag14}
\Pi(t,s;\gamma)\=\P(\mbox{no jump of } Y\mbox{ on }[t,t+s]|\gamma),\\
\label{eq:ag14a}
R(t,s;\gamma)\=\P(Y(t)=Y(t+s)|\gamma)
\end{eqnarray}
are more widely considered in the literature than their barred versions~(\ref{eq:ag2a}) and~(\ref{eq:ag6}) above.
In the present case there is almost surely no difference, since a.s.~$\gamma_x\ne\gamma_y$ provided $x\ne y$ and $\gamma_x>0$.
\end{rmk}

The above examples could be done either by considering the clock processes $\hat \G^{(\eps)}$ and $\hat\G$
on their own, together with~(\ref{eq:ag10}), in the case of Example 1, or, in the case of Example 2, we used,
besides Theorem~\ref{teo:asz}, convergence results for $S$ and $\G$ (in the appropriate scale), and in both examples
the limit is a correlation function of the limiting clock process $\hat\G$. Our last example 
is natural from the aging point of view, requires Theorem~\ref{teo:asz} alone, and the limit is not a function of $\hat\G$ alone. 

\paragraph{Example 3} 
Let
\begin{equation}
\label{eq:ag16}
\textstyle Q(t,s;\gamma)=\P\left(\sup_{r\in[0,t]}Z(r)<\sup_{r\in[0,t+s]}Z(r)|\gamma\right).
\end{equation}
This function was suggested in~\cite{FIN} as a ``measure of the prospects for novelty in the system``.
$\hat Z$ is almost surely continuous in single deterministic times, so we have that
\begin{equation}
\label{eq:ag17}
\lim_{\eps\to0}Q(\eps t,\eps s;\gamma)=\P\left(\textstyle\sup_{r\in[0,t]}\hat Z(r)<\sup_{r\in[0,t+s]}\hat Z(r)\right)=:\hat Q(t,s),
\end{equation}           
since the function $1\{\sup_{r\in[0,t]}\Phi(r)<\sup_{r\in[0,t+s]}\Phi(r)\}$ is %almost surely 
continuous in $\Phi\in D$ for almost every $\Phi$ under the distribution of $\hat Z$. 
We note that $\hat Q(t,s)$ is a function of the ratio $t/s$ only, by the self similarity of $\hat Z$,
but an explicit expression is not available, as far as we know, as it is for $\hat\Pi(t,s)$.

\vspace{.5cm}

\noindent{\bf Acknowledgements}

LRF would like to thank the CMI, Universit\'e de Provence, Aix-Marseille I for hospitality and support during several visits in the last few years
where this and related projects were developed. VG thanks NUMEC-USP for hospitality.

\end{document}